\documentclass[smallcondensed,runningheads]{svjour3}
\usepackage{amsmath,amssymb,mathtools,listings}
\usepackage[all,cmtip]{xy}
\usepackage[dvips]{graphicx}
\journalname{\textbf{\textsf{Review article}}}
\newcommand*{\id}{{\mathrm{id}}}
\newcommand*{\Ab}{{\mathbb A}}

\newcommand*{\Lc}{{\mathcal L}}

\newcommand*{\R}{{\mathbb R}}

\newcommand*{\Mcc}{{\mathcal M}}

\newcommand*{\Eb}{{\mathsf E}}

\newcommand*{\rd}{{\mathrm d}}
\newcommand*{\strat}{\,{\stackrel{\text{\tiny \textsf{o}}}{\phantom{\_}}}\,\,}

\spnewtheorem{algorithm}{Algorithm}{\sf}{\rm}
\smartqed

\begin{document}

\title{\textsf{An introduction to SDE simulation}}
\author{\textsf{Simon J.A. Malham} \and \textsf{Anke Wiese}}
\authorrunning{\textsf{Malham and Wiese}}
\institute{Simon J.A. Malham \and Anke Wiese \at
Maxwell Institute for Mathematical Sciences \\
and School of Mathematical and Computer Sciences\\
Heriot-Watt University, Edinburgh EH14 4AS, UK \\
Tel.: +44-131-4513200\\
Fax: +44-131-4513249\\
\email{S.J.Malham@ma.hw.ac.uk}\\
\email{A.Wiese@ma.hw.ac.uk}}

\date{5th April 2010}

\voffset=10ex 

\maketitle
\begin{abstract}
We outline the basic ideas and techniques underpinning the simulation
of stochastic differential equations. In particular we focus on strong simulation
and its context. We also provide illustratory examples and sample matlab algorithms
for the reader to use and follow. Our target audience is advanced undergraduate and 
graduate students interested in learning about simulating stochastic differential equations.
We try to address the FAQs we have encountered. 
\keywords{stochastic simulation}
\subclass{60H10 \and 60H35}
\end{abstract}

\lstset{language=Matlab,basicstyle=\ttfamily}

\section{\textsf{Introduction}}\label{intro}

\subsection{\textbf{\textsf{When is a model stochastic?}}}
Often in modelling, we need to incorporate a phenomenon or influence that
seems random, whose behaviour we can only recognize as
statistically Gaussian. The prototypical example is behaviour
of molecular origin---the Brownian motion of a pollen particle
on a water surface. However the influences can be large scale,
for example a turbulent wind or atmospheric flow, or
thousands of people buying and selling millions of shares.

Imagine tracking and adjusting the flight of a
rocket after lift-off. If the rocket is buffeted by a random turbulent
wind, you might sensibly equip the rocket with stabilizers that kick-in if a gust
diverts it too far. Computing the path of the rocket, and
regulating it to ensure it threads pre-arranged target positions
at critical junctures (eg.\/ stage separation), is
a stochastic simulation problem. Indeed it is a 
\emph{strong simulation problem}, as conditional on the path, 
the stabilizers will affect the outcome.
The pricing of financial derivatives (futures/options)
is another example. One tracks a security price (eg.\/ a share price) 
that is randomly buffeted by market forces. 
Pricing the derivative you wish to sell,
which might be exercised by the buyer at a future time,
involves hedging/regulating the proportion of your investment in
the security (the rest invested in risk-free bonds) so
as to minimize your risk. Building in stabilizers/barriers
that kick-in if the security price skews too wildly, 
is again, a strong stochastic simulation problem.

\subsection{\textbf{\textsf{What is a stochastic differential equation?}}}
Consider a model that incorporates some random phenomena whose
statistics is Gaussian. Suppose the state of the system is recorded through
the vector $y_t\in\R^N$, with $N\geqslant 1$. Suppose there are several random sources, say
$W^1,\ldots,W^d$; these are Wiener processes; think of them as independent, continuous,
nowhere differentiable, functions of time. Indeed the time derivatives of the Wiener processes
represent pure white noise. Suppose the affect of the Wiener processes on the model, 
i.e.\/ which way they skew the solution, is recorded through the vector fields
$V_1,\ldots,V_d$. Without the noise we would have a nice fuzz-free signal
which is generated by the vector field $V_0$. 
A stochastic model for the system state $y_t\in\R^N$ might 
be that it evolves according to the stochastic differential equation:
\begin{equation*}
\frac{\rd y}{\rd t}=V_0(y_t)
+V_1(y_t)\,\frac{\mathrm{d}W_t^1}{\rd t}+\cdots+V_d(y_t)\,\frac{\mathrm{d}W_t^d}{\rd t}.
\end{equation*}
This representation of the model is somewhat formal; after all the 
pure white noise terms $\rd W^i_t/\rd t$ need to interpreted in an extremely 
weak sense, we prefer to represent the model in the form
\begin{equation*}
\mathrm{d} y_t=V_0(y_t)\,\mathrm{d}t
+V_1(y_t)\,\mathrm{d}W_t^1+\cdots+V_d(y_t)\,\mathrm{d}W_t^d.
\end{equation*}
Indeed there is also a representation in a preferred integral form which
we meet presently. With this in mind though, an important observation at this point, 
is to recall that Wiener processes are continuous functions. Thus the  
solution function will also be continuous.\medskip

\noindent\textsf{Example (Langevin equation/Brownian motion).}
Consider the equation of motion of a pollen particle suspended
in a fluid flow. The particle might obey the following equation
of motion for its velocity $y_t$:
\begin{equation*}
\frac{\rd y_t}{\rd t}=-a\,y_t+\sqrt{b}\,\frac{\rd W_t}{\rd t},
\end{equation*}
where $a$ and $b$ are constants. The right-hand side is the force
exerted on the particle per unit mass. There is a deterministic
force $-a\,y_t$ and a white noise force $\sqrt{b}\,\rd W_t/\rd t$
supposed to represent the random buffeting by the water molecules.
\medskip

%

\subsection{\textbf{\textsf{What information do we want to retrieve?}}}
If the time interval of interest is $[0,T]$, and our initial deterministic state
is $y_0$, each realization $\omega$ of an individual Wiener path $W(\omega)$ 
will produce a different
outcome $y_t(\omega)$ for $t\in[0,T]$. Practical information of interest 
is often expected values of functions $f$ of the solution, $f(y_t)$,
or more generally path-dependent functions of the solution $f(t,y_s;s\leqslant t)$.
Hence we might want to compute
\begin{equation*}
\Eb\, f(y_t)\coloneqq\int f\bigl(y_t(\omega)\bigr)\,\rd\mathsf P(\omega),
\end{equation*}
where $\mathsf P$ is a probability measure.
For example we could pick $f$ to be of polynomial or exponential form, synomynous
with statistical moments of $y_t$. If $f$ is the identity map, we obtain
the expectation of the solution. If we take $f$ to be $\|\cdot\|_p^p$,
where $\|\cdot\|_p$ is the $p$-vector norm, then we define the $L^p$-norm for
$p\geqslant 1$ by
\begin{equation*}
\|y_t\|_{L^p}^p\coloneqq\int \|y_t(\omega)\|_p^p\,\rd\mathsf P(\omega).
\end{equation*}

\subsection{\textbf{\textsf{How do we retrieve it?}}}
There are two main simulation approaches to extract such information, we can either:
\begin{itemize}
\item Solve a \emph{partial differential equation}; or 
\item Perform \emph{Monte--Carlo simulation}.
\end{itemize}

Associated with every stochastic differential equation, there is a parabolic \emph{partial differential
equation} for $u(t,y)$ whose solution at time $t\in[0,T]$ is 
\begin{equation*}
u(t,y)=\Eb\, f(y_t)
\end{equation*}
provided $u(0,y)=f(y)$ initially.  Thus solving the associated partial differential equation
on $[0,T]$ will generate lots of information about the solution to the stochastic differential
equation at time $t$. By dudiciously choosing $f$ to be a monomial function we can generate
any individual moment of the solution $y_t$ we like, or if we choose $f=\exp$ we generate
all the moments simultaneously (this is essentially the Laplace transform). If we choose $f$ to be
a Dirac delta function we generate the transition probability distribution for $y_t$---the 
probability density function for $y_t$ conditioned on the initial data $y_0$. Choosing $f$
to be a Heaviside function generates the corresponding (cumulative) distribution function.
Of course often, the partial differential equation will have to be solved approximately. 
Also note, if we fix a form for $f$ from the start, for example $f$ as the identity map, 
then we simply solve an ordinary differential equation for $u(t,y)$.

In \emph{Monte--Carlo simulation}, we generate a set of suitable 
multidimensional sample paths 
$\hat W(\omega)\coloneqq\bigl(\hat W^1(\omega),\ldots,\hat W^d(\omega)\bigr)$ on $[0,T]$; 
in practice, $\omega$ belongs to a large but finite set.
For each sample path $\hat W(\omega)$, we generate a sample path solution $\hat y(\omega)$
to the stochastic differential equation on $[0,T]$.
This is often achieved using a truncation of the `stochastic' Taylor
series expansion for the solution $y$ of the stochastic differential
equation, on successive small subintervals of $[0,T]$.
Suppose for example, we wanted to compute the expectation $\Eb\,f(\hat y_t)$.
Having generated a set of approximate solutions 
$\hat y_t(\omega_i)$ at time $t\in[0,T]$, for $i=1,\ldots,P$ with $P$ large,
we can estimate $\Eb\,f(\hat y_t)$ by computing the mean-sum over the large finite 
set of approximate sample solutions $\hat y_t(\omega_i)$. 
Hence in practice we approximate
\begin{equation*}
\int f\bigl(y_t(\omega)\bigr)\,\rd\mathsf P(\omega)
\approx\tfrac{1}{P}\sum_{i=1}^Pf\bigl(y_t(\omega_i)\bigr)
\end{equation*}
where $P$ is the total number of sample paths.
A natural dichotomy now arises. 
To compute $\Eb\,f(\hat y_t)$, we can in fact choose any suitable
multidimensional paths $\hat W(\omega)$ that leave $\Eb\,f(\hat y_t)$
approximately invariant, in the sense that $\|\Eb\,f(y_t)-\Eb\,f(\hat y_t)\|$
is sufficiently small. This is a \emph{weak approximation}. 
For example, increments $\Delta W^i$ in each computation interval
can be chosen from a suitable binomial branching process, or using
Lyons and Victoir's cubature method~\cite{LV}. Note that since the 
approximate paths are not close to Brownian paths we cannot compare
$\hat y_t$ and $y_t$ directly. 
In a \emph{strong approximation}, discrete increments $\Delta W^i$ 
in each computation interval are directly sampled from the Gaussian
distribution. This is more expensive. However, the sample paths 
$\hat W(\omega)$ generated in this way, allow us to
compare $\hat y_t(\omega)$ and $y_t(\omega)$ directly in the sense that
we can guarantee $\Eb\,\|y_t-\hat y_t\|$ will be sufficiently small.
Naturally, using strong simulation we can also
account for path-dependent features, such
as conditional cut-offs or barriers, when we investigate
individual solutions or the final expectation or higher
moments of the approximate solution paths $\hat y$. 
For a comprehensive overview of Monte--Carlo methods
see Boyle, Broadie and Glasserman~\cite{BBG}.

\subsection{\textbf{\textsf{What is required?}}}
In general to extract qualitative and quantative information from a stochastic
differental system requires the languages and techniques of several mathematical
disciplines, notably:
\begin{enumerate}
\item \emph{Integration}: in Brownian motion new information is continuously 
generated on infinitesimally small time scales (imagine the pollen particle jiggles);
solution as with ordinary differential equations is by integration, except that 
now the coefficients of the evolution equation---the Wiener processes---are no longer differentiable.
\item \emph{Statistics}: we typically extract statistical information from the solution
process;
\item \emph{Geometry:} as with ordinary differential equations, preserving invariant 
geometric structure of the solution path evolution is important; for example the solution
may evolve on a homogeneous manifold;
\item \emph{Simulation:} stochastic differential equations more often than not,
are \emph{not} integrable in the classical sense, and require numerical computation. 
\end{enumerate}
For general background reading, we recommend as follows. For a comprehensive 
introduction to the theory underlying stochastic differential equations download
Evans' notes~\cite{Evans}. For an introduction to numerical simulation, see
Higham's notes~\cite{Higham:intro}.
The answer to just about any other question that a beginner may have on numerical simulation,
not covered above or here, can likely be found in the treatise by 
Kloeden and Platen~\cite{KP}.

\section{\textsf{Stochastic differential equations}}

\subsection{\textbf{\textsf{Integral representation}}}
Consider the nonlinear stochastic differential equation of order $N\in\mathbb N$ given by 
\begin{equation*}
y_t=y_0+\int_0^t \tilde V_0(y_\tau)\,\mathrm{d}\tau
+\sum_{i=1}^d\int_0^t V_i(y_\tau)\,\mathrm{d}W^i_\tau.
\end{equation*}
Here $(W^1,\ldots,W^d)$ is a $d$-dimensional Wiener process, i.e.\/
there are $d$ independent driving noisy signals.
We assume there exists a unique solution $y\colon[0,T]\mapsto\mathbb R^N$ 
for some time interval $[0,T]\subseteq\mathbb R_{+}$. 
We suppose that $\tilde V_0$ and $V_i\colon\mathbb R^N\rightarrow \mathbb R^N$, 
$i=1,\ldots,d$, are smooth non-commuting autonomous vector fields.
We are representing the stochastic differential equation above in It\^o form
and indicate this by using $\tilde V_0$ to represent the \emph{It\^o drift vector field}.
We call the vector fields $V_i$ for $i=1,\ldots,d$ associated with the
driving noise terms the \emph{diffusion vector fields}.
Presently we will distinguish, and explain, the It\^o
representation as opposed to the Stratonovich representation for a stochastic
differential equation. We also remark that a common convention is to set $W^0_t\equiv t$.
Results on existence and uniquess of solutions can be found in Kloeden and Platen~\cite{KP}.

\subsection{\textbf{\textsf{Driving Wiener process}}}
A scalar driving noisy signal or disturbing Brownian motion
has a concise definition and set of properties formulated by Wiener.

\begin{definition}[\textbf{\textsf{Wiener process}}]
A scalar \emph{standard Wiener process} or \emph{standard Brownian motion} $W$ 
is a continuous process that satisfies the three conditions:
\begin{enumerate}
\item $W_0=0$ with probability one;
\item $W_t-W_s\sim \sqrt{t-s}\,\cdot\text{\textsf{N}}(0,1)$
for $0\leqslant s<t$, where $\text{\textsf{N}}(0,1)$
denotes a standard Normal random variable;
\item Increments $W_t-W_s$ and $W_\xi-W_\eta$ on distinct time intervals 
are independent, i.e.\/ for $0\leqslant s<t<\eta<\xi$.
\end{enumerate}
Note that with probability one an individual Brownian path 
is nowhere differentiable.
\end{definition}
\medskip

\begin{figure}
  \begin{center}
  \includegraphics[width=7.0cm,height=5.0cm]{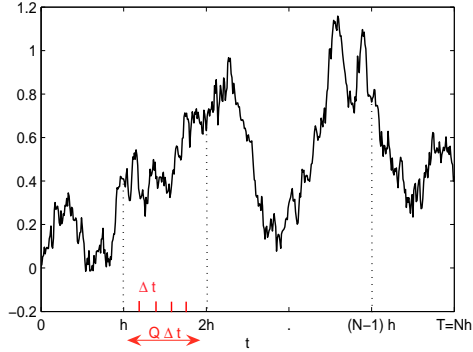}
  \end{center}
\caption{Example scalar Wiener path on the interval $[0,T]$.}
\end{figure}

\noindent\textsf{Example (Langevin equation).} As we have seen,
the Brownian motion of a pollen particle suspended
in a fluid flow obeys the following equation
of motion for its velocity $y_t$:
\begin{equation*}
\rd y_t=-a\,y_t\,\rd t+\sqrt{b}\,\rd W_t,
\end{equation*}
where $a$ and $b$ are constants, and $W$ is a scalar Wiener process.
This type of stochastic differential equation is said to have \emph{additive noise}
as the diffusion vector field is constant. It is also an example of an 
\emph{Ornstein--Uhlenbeck process}.
\medskip

\noindent\textsf{Example (scalar linear equation).} 
Consider the scalar linear stochastic differential equation
\begin{equation*}
\rd y_t=a\,y_t\,\rd t+b\,y_t\,\rd W_t
\end{equation*}
driven by a scalar Wiener process $W$, with $a$ and $b$ constants.
This stochastic differential equation is said to have \emph{multiplicative noise}
as the diffusion vector field depends multiplicatively on the solution $y_t$. We
can in fact analytically solve this equation, 
the solution is
\begin{equation*}
y_t=y_0\,\exp\bigl(a\,t+b\,W_t-\tfrac12b^2t\bigr).
\end{equation*}
The additional term `$-\tfrac12b^2t$' is due to the It\^o correction, which
we discuss shortly.\medskip

\subsection{\textbf{\textsf{Vector fields and flow maps}}}
Consider an ordinary differential equation governed by an autonomous vector field $V$
that evolves on a homogeneous manifold $\Mcc$ so that
\begin{equation*}
\frac{\mathrm{d}y_t}{\mathrm{d}t}=V(y_t),
\end{equation*}
with initial data $y_0\in\Mcc$. Here we suppose $\Mcc$ to be an 
$N$-dimensional manifold. The reader can for simplicity assume $\Mcc\equiv\R^N$
for the rest of this section if they choose. Let $\text{Diff}(\Mcc)$ denote the
group of diffeomorphisms of $\Mcc$. The \emph{flow-map} $\varphi_{t,t_0}\in\text{Diff}(\Mcc)$
for the ordinary differential equation above is the map taking the 
solution configuration on $\Mcc$ at time $t_0$ to that at time $t$, 
i.e.\/ it is the map $\varphi_{t,t_0}\colon\Mcc\to\Mcc$ such that
\begin{equation*}
\varphi_{t,t_0}\colon y_{t_0}\mapsto y_t.
\end{equation*}
In other words for any data $y_0\in\Mcc$ at time $t_0$ we can 
determine its value $y_t\in\Mcc$ at time $t$ later by applying the \emph{action}
of the flow map $\varphi_{t,t_0}$ to $y_0$ so that $y_t=\varphi_{t,t_0}\circ y_0$.
Note that the flow map satisfies the usual group properties
\begin{equation*}
\varphi_{t,s}\circ\varphi_{s,t_0}=\varphi_{t,t_0},
\end{equation*}
with $\varphi_{t_0,t_0}=\id$, the identity diffeomorphism.
If $f\in\text{Diff}(\Mcc)$, 
the chain rule reveals that for all $y\in\Mcc$, we have
\begin{equation*}
\frac{\mathrm{d}f(y)}{\mathrm{d}t}=V(y)\cdot\partial_y\,f(y),
\end{equation*}
where $\partial_y\equiv\nabla_y$ is the usual gradient operator with respect 
to each component of $y$. In other words, vector fields act on the group 
of diffeomorphisms $\text{Diff}(\Mcc)$ 
as first order partial differential operators, 
and for any $f\in\text{Diff}(\Mcc)$, we write
\begin{equation*}
V\circ f\circ y=V(y)\cdot\partial_yf(y).
\end{equation*}
We now think of $V(y)\cdot\partial_y$ as a first order partial differential operator
and an element of the tangent space $\mathfrak X(\Mcc)$ to $\text{Diff}(\Mcc)$.
In particular, choose the diffeomorphic map $f$ to be the flow map 
$\varphi_t\equiv\varphi_{t,0}$. Then for all $y_0\in\Mcc$ 
and $y_t=\varphi_t\circ y_0$ we have, 
\begin{equation*}
\frac{\rd}{\rd t}(\varphi_t\circ y_0)=V\circ\varphi_t\circ y_0.
\end{equation*}
We pull back this ordinary differential equation for $y_t\in\Mcc$ 
to the \emph{linear} functional differential equation in $\text{Diff}(\Mcc)$:
\begin{equation*}
\frac{\rd\varphi_t}{\rd t}=V\circ\varphi_t.
\end{equation*}
Since $\varphi_0=\id$, the solution is
$\varphi_t=\exp(t\,V)$, giving the represention of the flow-map as the 
\emph{exponentiation of the vector field}. Hence we see that
\begin{equation*}
y_t=\exp(t\,V)\circ y_0.
\end{equation*}

An important and illustrative concomitant derivation of this result is as follows.
Integrating the functional differential equation for the flow-map we get
\begin{equation*}
\varphi_t=\id+\int_0^t V\circ\varphi_\tau\,\rd\tau.
\end{equation*}
To solve this integral equation, we set up the formal iterative procedure given by
\begin{equation*}
\varphi_t^{(n+1)}=\id+\int_0^t V\circ\varphi^{(n)}_\tau\,\rd\tau,
\end{equation*}
with $\varphi^{(0)}_t=\id$.
For example, after two iterations:
$\varphi_t^{(2)}=\id+t\,V\circ\id+\tfrac12t^2\,V^2\circ\id$,
where $V^2\equiv V\circ V$. Hence in the limit we obtain the exponential form
for $\varphi_t$ above.

The composition of two vector fields $U$ and $V$ is a second order differential
operator:
\begin{align*}
U\circ V=&\;\bigl(U(y)\cdot\partial_y\bigr)\bigl(V(y)\cdot\partial_y\bigr)\\
=&\;\Bigl(\bigl(U(y)\cdot\partial_y\bigr)\bigl(V(y)\bigr)\Bigr)\cdot\partial_y
+U(y)\otimes V(y)\colon \partial_{yy}\\
=&\;\sum_{i,j=1}^N U^i\partial_{y_i}(V^j)\partial_{y_j}+
\sum_{i,j=1}^N U^iV^j\partial_{y_iy_j}.
\end{align*}
Importantly we now observe that since $\partial_{y_iy_j}=\partial_{y_jy_i}$ 
as operators on $\text{Diff}(\Mcc)$, 
the Lie bracket $[U,V]$ of two vector fields is a vector field:
\begin{equation*}
[U,V]=U\circ V-V\circ U
=\Bigl(\bigl(U(y)\cdot\partial_y\bigr)V(y)-
\bigl(V(y)\cdot\partial_y\bigr)U(y)\Bigr)\cdot\partial_y.
\end{equation*}
Note the sum of two vector fields is itself a vector field. Hence
the set of vector fields 
is a \emph{Lie algebra} $\mathfrak X(\Mcc)$---closed under summation
and Lie product $[\cdot,\cdot]$.

\subsection{\textbf{\textsf{Stratonovich representation}}}
There are two generic representations for stochastic differential
equations. One can either express them in It\^o form, as we did at
the beginning of Section~2, or we can express the stochastic differential
equation in Stratonovich form, in which case we write
\begin{equation*}
y_t=y_0+\int_0^t V_0(y_\tau)\,\mathrm{d}\tau
+\sum_{i=1}^d\int_0^t V_i(y_\tau)\,\strat\mathrm{d}W^i_\tau.
\end{equation*}
Two subtle notational changes can be spotted. First the 
`$\strat\mathrm{d}W^i_\tau$' indicates that the
stochastic integrals on the right are supposed to be interpreted
in the Stratonovich sense; discussed presently. Second we
now use the \emph{Stratonovich drift vector field} $V_0$ instead
of the It\^o drift vector field $\tilde V_0$. 
The relation between the two vector fields is
\begin{equation*}
V_0=\tilde V_0-\tfrac12\sum_{i=1}^d(V_i\cdot\partial_yV_i).
\end{equation*}
Importantly, when stochastic integrals are interpreted in the It\^o sense,
then they are limit of a left Riemann sum and when repeated integrals are
computed, an It\^o correction must be taken into account; a practical
discussion of this point can be found in Higham~\cite[pp.~530--1]{Higham:intro}.  
For example, the correct evaluation of the an It\^o repeated integral 
of a Wiener process with respect to itself is
\begin{equation*}
\int_0^T W^i_\tau\,\mathrm{d}W^i_\tau=\tfrac12\bigl(W^i_T\bigr)^2-\tfrac12 T.
\end{equation*}
Stratonovich integrals are interpreted as the limit of the midpoint rule,
so for the corresponding Stratonovich integral the correct evaluation is
\begin{equation*}
\int_0^T W^i_\tau\,\strat\mathrm{d}W^i_\tau=\tfrac12\bigl(W^i_T\bigr)^2.
\end{equation*}
Thus the rules of Stratonovich integral calculus match those of 
standard integral calculus. For this reason it is often preferable 
to use the Stratonovich representation for a stochastic differential equation.
The two representations are equivalent, but it is important to know
in which form you have been quoted the stochastic differential equation.
Often this depends on the modeller and their field of interest.
In finance applications the It\^o representation 
predominates, and by simply replacing the given It\^o drift vector field $\tilde V_0$
by the \emph{corresponding} Stratonovich drift vector field $V_0$ above, one can
proceed using standard integral calculus rules. In physical applications,
the model is often often directly expressed in Stratonovich form.
Hereafter we will use the Stratonovich representation and omit the `$\strat$' symbol, 
unless we specify otherwise.

\section{\textsf{Stochastic Taylor expansion}}
We follow the procedure we performed above for the ordinary differential equation 
to try to find the solution for the flow-map. In the process
we obtain a solution series expansion called the stochastic Taylor series.
We define the \emph{flow-map} $\varphi_{t}\in\text{Diff}(\Mcc)$
for the stochastic differential equation above as the map taking the 
solution configuration on $\Mcc$ at time $0$ to that at time $t$; hence
$y_t=\varphi_t\circ y_0$.
Using the Stratonovich representation for a stochastic differential
equation and the convention $W^0_t\equiv t$, 
the chain rule for any function $f\in\text{Diff}(\Mcc)$ yields the 
stochastic differential equation governing the evolution
of $f\circ y_t$ as follows
\begin{equation*}
f\circ y_t=f\circ y_0
+\sum_{i=0}^d\int_0^t (V_i\cdot\partial_y\,f)\circ y_\tau\,\mathrm{d}W^i_\tau.
\end{equation*}
As for the ordinary differential equation, setting $f=\varphi_t$, 
we can pull back the stochastic
flow on $\Mcc$ to a functional stochastic differential equation 
on $\text{Diff}(\Mcc)$ given by
\begin{equation*}
\varphi_t=\id+\sum_{i=0}^d\int_0^t V_i\circ\varphi_\tau\,\mathrm{d}W^i_\tau.
\end{equation*}
To solve this equation, we set up the formal iterative procedure given by
\begin{equation*}
\varphi_t^{(n+1)}=\id+\sum_{i=0}^d\int_0^t V_i\circ\varphi^{(n)}_\tau\,\mathrm{d}W^i_\tau,
\end{equation*}
with $\varphi^{(0)}_t=\id$.
By performing the iterations one can see formally, and prove rigorously,
that the solution flow-map is given by the series expansion
\begin{equation*}
\varphi_t=\id+\sum_{i=0}^d\bigl(W_t^i\bigr)\,V_i
+\sum_{i,j=0}^d\biggl(\int_0^t\int_0^{\tau_1}\,\rd W^i_{\tau_2}\,\rd W^j_{\tau_1}\biggr)\,V_{ij}+\cdots.
\end{equation*}
Here we use the notation $V_{ij}\equiv V_i\circ V_j$. We can apply this to the 
initial data $y_0\in\Mcc$ and obtain the \emph{stochastic Taylor expansion} 
for the solution
\begin{equation*}
y_t=y_0+\sum_{i=0}^d\bigl(W_t^i\bigr)\,V_i(y_0)
+\sum_{i,j=0}^d\biggl(\int_0^t\int_0^{\tau_1}\,\rd W^i_{\tau_2}\,\rd W^j_{\tau_1}\biggr)\,V_{ij}(y_0)+\cdots.
\end{equation*}
We can express the solution series for the flow-map concisely 
as follows. Let $\Ab^\ast$ denote the free monoid of words over the
alphabet $\mathbb A=\{0,1,\ldots,d\}$. We adopt the standard notation for 
\emph{Stratonovich integrals}, if $w=a_1\ldots a_n$ then we set
\begin{equation*}
J_w(t)\coloneqq\int_0^t\cdots\int_0^{\tau_{n-1}}
\mathrm{d}W^{a_1}_{\tau_n}\,\cdots\,\mathrm{d}W^{a_n}_{\tau_1}.
\end{equation*}
We also write the composition of the vector fields as
$V_w\equiv V_{a_1}\circ V_{a_2}\circ\cdots\circ V_{a_n}$. Then the flow-map
is given by
\begin{equation*}
\varphi_t=\sum_{w\in\Ab^\ast}J_w(t)\,V_w.
\end{equation*}

\section{\textsf{PDE simulation}}
There is an intimate link between any stochastic differential
equation and a prescribed \emph{parabolic partial differential equation}.
The link is given by the Feynman--Kac formula, which we give here
in a very simple form. See for example Karlin and Taylor~\cite[pp.~222--4]{KT} 
for the full statement of the Feynman--Kac formula and its applications.

\begin{theorem}[\textbf{\textsf{Feynman--Kac formula}}]
Consider the parabolic partial differential equation for $t\in[0,T]$:
\begin{equation*}
\partial_t u=\Lc\, u,
\end{equation*}
with $u(0,y)=f(y)$. Here 
$\mathcal L\coloneqq V_0+\tfrac12(V_1^2+\cdots+V_d^2)$ 
is a differential operator of order $2N$.
Let $y_t$ denote the solution
to the stochastic differential equation for $t\in[0,T]$:
\begin{equation*}
y_t=y_0+\int_0^t V_0(y_\tau)\,\mathrm{d}\tau
+\sum_{i=1}^d\int_0^t V_i(y_\tau)\,\mathrm{d}W^i_\tau.
\end{equation*}
Then, when $y_0=y$ we have: $u(t,y)=\Eb\,f(y_t)$.
\end{theorem}
\medskip

\noindent\textsf{Remark.}
Note that using the relation between the It\^o and Stratonovich drift vector fields, 
an equivalent formulation is
$\Lc\equiv \tilde V_0\cdot\partial_y+\tfrac12\sum_{i=1}^d(V_i\otimes V_i):\partial_{yy}$.
\medskip

We provide a purely combinatorial proof of the Feynman--Kac formula.
Before we begin, we need the following results, for the expectation of 
Stratonovich integrals and also a combinatorial expansion. 
Let $\mathbb D^*\subset\mathbb A^*$ denote the 
free monoid of words constructed from the alphabet $\mathbb D=\{0,11,22,\ldots,dd\}$.
The expectation of a Stratonovich integral $J_w$ is given by
\begin{equation*}
\Eb\,J_w=\begin{cases} 
 \frac{t^{\mathrm{n}(w)}}{2^{\mathrm{d}(w)}\mathrm{n}(w)!},&\qquad w\in\mathbb D^*; \\
 0,&\qquad w\in\Ab^\ast\backslash\mathbb D^*.
\end{cases}
\end{equation*}
In the formula, $\mathrm{d}(w)$ is the number of non-zero consecutive pairs from $\mathbb D$
in $w$ and $\mathrm{n}(w)=\mathrm{z}(w)+\mathrm{d}(w)$, where $\mathrm{z}(w)$ is
the number of zeros in $w$. 
We also have the following combinatorial identity for all $w\in\mathbb D^\ast$:
\begin{equation*}
\bigl(V_0+\tfrac12(V_1^2+\cdots+V_d^2)\bigr)^k
\equiv\sum_{\mathrm{n}(w)=k}(\tfrac12)^{\mathrm{d}(w)}\,V_w,
\end{equation*}
where note that $V_i^2\equiv V_{ii}$.
In other words, expanding $\bigl(V_0+\tfrac12(V_{1}^2+\cdots+V_{d}^2)\bigr)^k$ generates
all the possible vector fields $V_w$ with $w\in\mathbb D^*$ and $\mathrm{n}(w)=k$, with
the appropriate coefficients of powers of one-half.

\begin{proof}
In the series solution for the flow-map $\varphi_t$, 
all stochastic information is encoded 
in the words on the left and the geometric information on the right.
Taking the expectation of the flow-map, 
noting that expectation is a linear operator,
and using the two results above for $\Eb\,J_w$ and the combinatorial expansion,
we get
\begin{align*}
\Eb\,\varphi_t=&\;\Eb\,\sum_{w\in\Ab^*} J_wV_w\\
=&\;\sum_{w\in\Ab^*}\bigl(\Eb\,J_w\bigr)\,V_w,\\
=&\;\sum_{k\geq0} \sum_{\mathrm{n}(w)=k}(\tfrac12)^{\mathrm{d}(w)}\frac{t^k}{k!}\,V_w\\
=&\;\sum_{k\geq0} \frac{t^k}{k!}\bigl(V_0+\tfrac12(V_{1}^2+\cdots+V_{d}^2)\bigr)^k\\
=&\;\exp\bigl(t\,\mathcal L\bigr).
\end{align*}
Now note that $\exp\bigl(t\,\mathcal L\bigr)$ generates the \emph{semi-group} 
for the solution to the parabolic differential equation in the theorem.\qed
\end{proof}
\medskip

\noindent\textsf{Example (Heston model).}
In the Heston model~\cite{Heston}, a stock price $S_t$ is modelled by 
a stochastic process $x_t=\log S_t$ with variance process $v_t$ 
which evolve according to:
\begin{align*}
\mathrm{d}x_t=&\;\mu\,\mathrm{d}t+\sqrt{v_t}\,\mathrm{d}W^1_t,\\
\mathrm{d}v_t=&\;\kappa(\theta-v_t)\,\mathrm{d}t
+\varepsilon\,\sqrt{v_t}\,\bigl(\rho\,\mathrm{d}W^1_t
+\sqrt{1-\rho^2}\,\mathrm{d}W^2_t\bigr),
\end{align*} 
given in It\^o form.
Note that the variance is a mean-reverting process; it tries to revert
to the mean value $\theta$ at rate $\kappa$. Using the Feynman--Kac formula
the corresponding partial differential equation for 
$u(x,v,t)=\Eb\,\bigl(f(x_t,v_t)~|~x_0=x,~v_0=v\bigr)$ is
\begin{equation*}
u_t=\mu\,u_x+\kappa(\theta-v)\,u_v
+\tfrac12v\,u_{xx}+\rho\epsilon v\,u_{xv}+\tfrac12\epsilon^2v\,u_{vv}.
\end{equation*}
\medskip

\noindent\textsf{Remark.} The Feynman--Kac formula shows how
we can solve a partial differential equation to obtain information 
about the solution $y_t$ to the stochastic differential equation at time $t\in[0,T]$. 
In the reverse direction, to numerically solve high dimensional diffusion problems
in the form of deterministic partial differential equations, we need only 
simulate an $N$-dimensional stochastic equation for $y_t$ and then compute 
the expectation $\Eb\,y_t$ to find the solution.\smallskip

\section{\textsf{Monte--Carlo simulation}}
In \emph{Monte--Carlo simulation}, we generate a set of suitable  
multidimensional sample paths say 
$\hat W(\omega)\coloneqq\bigl(\hat W^1(\omega),\ldots,\hat W^d(\omega)\bigr)$ 
on $[0,T]$. We generate a large finite set of paths, each labelled by 
$\omega$. Here the number of paths, say $P$,
must be large enough so that, for example, any statistical information 
for the solution $y_t$ that we want to extract is sufficiently robust. 
For each sample path $\hat W(\omega)$, we generate a sample path solution $\hat y(\omega)$
to the stochastic differential equation on $[0,T]$. This can often only
be achieved approximately, by using a truncation of the stochastic Taylor
expansion for the solution $y$ on successive small subintervals of $[0,T]$.
Having generated a set of approximate solutions 
$\hat y_t(\omega_i)$ at time $t\in[0,T]$ for every $\omega_i$ for $i=1,\ldots,P$, 
we estimate the expectation $\Eb\,f(\hat y_t)$ by computing
\begin{equation*}
\tfrac{1}{P}\sum_{i=1}^Pf\bigl(\hat y_t(\omega_i)\bigr)
\end{equation*}
regarded as a suitable approximation for 
\begin{equation*}
\Eb\,f(y_t)\coloneqq\int 
f\bigl(y_t(\omega)\bigr)\,\rd\mathsf P(\omega)
\end{equation*}
over all possible paths.
Now a natural question arises. Do the suitable paths $\hat W(\omega)$ we generate, 
have to be sample Brownian paths to compute the mean above, or
can we choose different paths that will still generate the expectation effectively?
We discuss the latter case (weak simulation) briefly next. 
We then move onto a more comprehensive treatment of the former 
(strong simulation) case, which also allows us to include
path-dependent features.

\subsection{\textbf{\textsf{Weak simulation}}}
There are several \emph{weak simulation} strategies to approximate $\Eb\,f(\hat y_t)$.
The simplest and most common is to replace the driving Wiener paths $W^i$
by paths generated as follows.
Construct paths by generating increments $\Delta W^i(t_n,t_{n+1})$
over the computation subinterval $[t_n,t_{n+1}]$ by the binomial branching process
\begin{equation*}
\mathsf{P}\bigl(\Delta W^i(t_n,t_{n+1})=\pm\sqrt{h}\bigr)=\tfrac12,
\end{equation*}
where $h=t_{n+1}-t_n$. 
Then, depending on the ordinary numerical integration scheme employed
for each such path, one can show that for some order of convergence $p$,
for some $t\in[0,T]$, we have 
\begin{equation*}
\|\Eb\,f(y_t)-\Eb\,f(\hat y_t)\|=\mathcal O(h^p).
\end{equation*}
See Kloeden and Platen~\cite{KP} for more details.
Another promising method is the cubature method of Lyons and Victoir~\cite{LV}.

\subsection{\textbf{\textsf{Strong simulation}}}
In a \emph{strong simulation}, discrete increments $\Delta W^i$ 
in each computation interval are directly sampled from the Gaussian
distribution. Indeed we generate each multidimensional
path $\hat W(\omega)$ by choosing increments in each component as follows
\begin{equation*}
\Delta W^i(t_n,t_{n+1})\sim\sqrt{h}\,\cdot\textsf{N}(0,1).
\end{equation*}
This is more expensive, but sample paths $\hat W(\omega)$ generated in this way, 
allow us to compare $\hat y_t(\omega)$ and $y_t(\omega)$ directly in the sense that
one can show
\begin{equation*}
\Eb\,\|y_t-\hat y_t\|=\mathcal O\bigl(h^{\frac{p}{2}}\bigr)
\end{equation*}
for some order of strong convergence $p/2$; which we will discuss in
detail in Section~\ref{sec:strongerror}. 
Often in practice we take $\|\cdot\|$ to be the
Euclidean norm so that the convergence shown is in the $L^2$-norm.

Given a sample mulitdimensional path $\hat W(\omega)$ on $[0,T]$,
how do we actually construct an approximate solution $\hat y_t$?
Here we are guided by the stochastic Taylor expansion. Indeed,
classical strong numerical methods are based on truncating the 
stochastic Taylor expansion 
\begin{equation*}
y_t=\sum_{w\in\Ab^\ast}J_w(t)\,V_w(y_0),
\end{equation*}
and applying the approximation over successive
subintervals of the global interval of integration $[0,T]$; 
see Kloeden and Platen~\cite{KP} or Milstein~\cite{Mil}.
We present three simple example numerical approximation methods.

\subsection{\textbf{\textsf{Euler--Maruyama method}}}
If we truncate the It\^o form of the stochastic Taylor series
after the first order terms we generate the 
\emph{Euler--Maruyama numerical method} as follows:
\begin{equation*}
\hat y_{n+1}=\hat y_n+h\,\tilde V_0(\hat y_n)
+\sum_{i=1}^d\bigl(\Delta W^i(t_n,t_{n+1})\bigr)\,V_i(\hat y_n).
\end{equation*}
This is a numerical scheme with global order of convergence
of $h^{\frac12}$. We explain in Section~\ref{sec:strongerror}
why we have used the It\^o drift vector field here.


\begin{figure}
  \begin{center}
  \includegraphics[width=12.0cm,height=9.0cm]{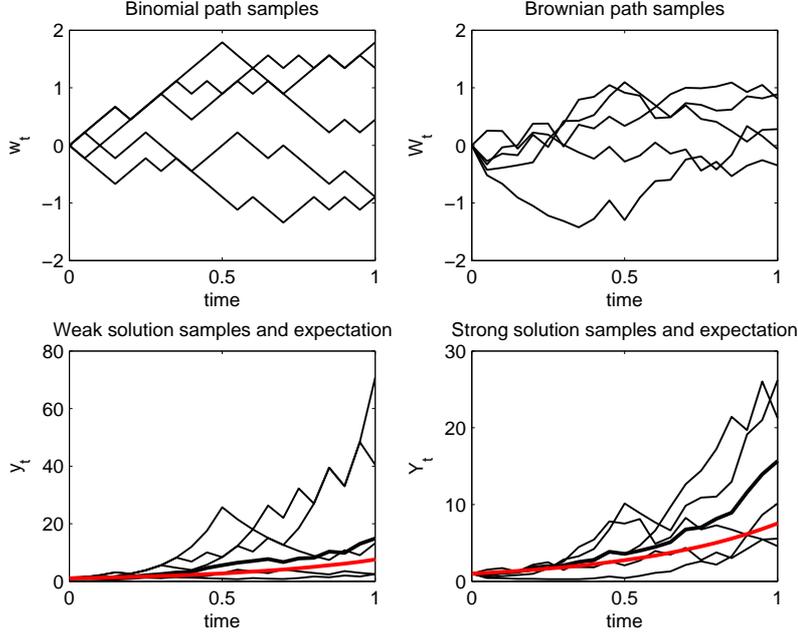}
  \end{center}
\caption{Weak and strong simulations for the scalar linear example
given in Section~2, with $a=3$ and $b=2$. For this example we took
$y_0=1$, $T=1$, $h=0.05$ and $P=10$ for both simulations, though only
$5$ sample paths are shown above in all cases. Top left are $5$ sample
binomial branching paths $w$, and top right are $5$ sample Brownian
paths $W$. Lower left are $5$ sample solution paths $y$ using the binomial branching
paths, while lower right are $5$ sample solution paths $Y$ 
using the Brownian paths; both computed 
using the Euler--Maruyama method. At each time-step the thick black line shows the
average value over $P=10$ samples and the red line is the analytic solution
for the expectation.}
\label{fig:weakstrongsim}
\end{figure}

\subsection{\textbf{\textsf{Milstein method}}}
We now truncate the stochastic Taylor series after the second order
terms. This generates the \emph{Milstein numerical method} given by
\begin{equation*}
\hat y_{n+1}=\hat y_n+h\,V_0(\hat y_n)+\sum_{i=1}^d\bigl(\Delta W^i(t_n,t_{n+1})\bigr)\,V_i(\hat y_n)
+\sum_{i,j=1}^dJ_{ij}(t_n,t_{n+1})V_{ij}(\hat y_n).
\end{equation*}
An important and expensive ingredient in this method,
is the simulation of the multiple integrals $J_{ij}(t_n,t_{n+1})$ 
for $i\neq j$ shown, on each integration step. When $i=j$,
the multiple integrals $J_{ii}(t_n,t_{n+1})$ are cheaply
evaluated by
\begin{equation*}
J_{ii}(t_n,t_{n+1})=\int_{t_n}^{t_{n+1}}\int_{t_n}^{\tau_1}\,\rd W^i_{\tau_2}\,\rd W^j_{\tau_1}=
\tfrac12\bigl(\Delta W^i(t_n,t_{n+1})\bigr)^2.
\end{equation*}
When $i\neq j$ we have by integration by parts that 
\begin{equation*}
J_{ji}=J_iJ_j-J_{ij}.
\end{equation*}
Hence we need only compute one double integral for each pair $i\neq j$. 
Equivalently we need only compute the L\'evy area given by
\begin{equation*}
A_{ij}\coloneqq\tfrac12(J_{ij}-J_{ji}),
\end{equation*}
since $J_{ij}=\tfrac12J_iJ_j+A_{ij}$. By Stokes' Theorem, 
the L\'evy area on the interval $[t_n,t_{n+1}]$ is the chordal area 
for the path $(W^i,W^j)$ on $[t_n,t_{n+1}]$. This can be observed
directly by definition since
\begin{equation*}
A_{ij}(t)=\tfrac12\int_{0}^{t}W^i_{\tau}\,\rd W^j_{\tau}-W^j_{\tau}\,\rd W^i_{\tau}.
\end{equation*}
We consider the issue of simulating the L\'evy area
in some detail in Section~\ref{sec:levy}. The Milstein scheme
has global order of convergence $h$.

\subsection{\textbf{\textsf{Castell--Gaines method}}}
Consider the exponential Lie series $\psi_t=\log\varphi_t$ 
generated by taking the logarithm of the 
stochastic Taylor series for the flow-map, i.e.\/
\begin{align*}
\psi_t&=(\varphi_t-\id)-\tfrac12(\varphi_t-\id)^2+\tfrac13(\varphi_t-\id)^3+\cdots\\
&=\sum_{i=0}^dJ_iV_i+\sum_{i>j}\tfrac12(J_{ij}-J_{ji})[V_i,V_j]+\cdots.
\end{align*}
This series is also known as the Chen--Strichartz, Chen--Fleiss or Magnus series.
The Castell--Gaines method is a strong numerical method 
based on truncating the exponential Lie series. As for the methods above, we
generate a set of multidimensional paths $\hat W(\omega)$ on $[0,T]$ with 
Wiener increments $\Delta W^i(t_n,t_{n+1})$ sampled on the scale $h=t_{n+1}-t_n$.
On each computation interval $[t_n,t_{n+1}]$, we relace the $J_i(t_n,t_{n+1})$
by the Normal samples $\Delta W^i(t_n,t_{n+1})$. If required, 
the L\'evy area increments $A_{ij}(t_n,t_{n+1})$ shown, 
are also replaced by suitable samples $\hat A_{ij}(t_n,t_{n+1})$ 
as we outline in Section~\ref{sec:levy}. Then across the 
computation interval $[t_n,t_{n+1}]$, we have 
\begin{equation*}
\hat\psi_{t_n,t_{n+1}}=\sum_{i=0}^d\bigl(\Delta W^i(t_n,t_{n+1})\bigr)\,V_i
+\sum_{i>j}\hat A_{ij}(t_n,t_{n+1})\,[V_i,V_j].
\end{equation*}
The solution at time $t_{n+1}$ is then approximately given by 
\begin{equation*}
\hat y_{t_{n+1}}\approx\exp(\hat\psi_{t_n,t_{n+1}})\circ\hat y_{t_n}.
\end{equation*}
Note that for each path $\Delta W^i(t_n,t_{n+1})$ and $\hat A_{ij}(t_n,t_{n+1})$ are 
fixed constants. Hence the truncated Lie series $\hat\psi_{t_n,t_{n+1}}$ 
is itself an autonomous vector field.
Thus, for $\tau\in[0,1]$ and with $u(0)=\hat y_{t_n}$,
we solve the ordinary differential equation
\begin{equation*}
u'(\tau)=\hat\psi_{t_n,t_{n+1}}\circ u(\tau).
\end{equation*}
Using a suitable high order ordinary
differential integrator generates 
$u(1)\approx\hat y_{t_{n+1}}$.

Without the L\'evy area the Castell--Gaines method has global
order of convergence $h^{\frac12}$; while with the L\'evy area it
has global order of convergence $h$.
Castell and Gaines~\cite{CG1,CG2} prove that
their strong order $h^{\frac12}$ method is
always more accurate than the Euler--Maruyama method. Indeed
they prove that this method is \emph{asymptotically efficient}
in the sense of Newton~\cite{Newton}. Further in the case of a single
driving Wiener process ($d=1$), they prove the same is true for
their strong order $h$ method. By asymptotically
efficient we mean, quoting from Newton, that they 
``minimize the leading coefficient in the expansion of mean-square errors
as power series in the sample step size''.

\section{\textsf{Simulating the L\'evy area}}\label{sec:levy}
A fundamental and crucial aspect to the implementation of
strong order one or higher integrators for stochastic differential 
equations, is the need to successfully simulate the L\'evy chordal
areas $A_{ij}(t_n,t_{n+1})$, when the diffusion vector fields do not commute.
This aspect is more than just a new additional concern once we step off the cliff edge
of simple path increment approximations with frozen vector fields characterized by
the Euler--Maruyama approximation. It also
represents a substantial technical difficulty. Here we will outline several
methods employed to simulate it sufficiently accurately; the important distinguishing 
criterion for the success of the method will be its asymptotic rate of convergence
as $h\to0$. A sample survey of these methods can be found in 
Ryd\'en and Wiktorsson~\cite{RW}. Here we will focus on the case of two
independent Wiener processes $W^1$ and $W^2$ and the requirement to 
simulate $A_{12}(h)\coloneqq A_{12}(t_n,t_{n+1})$, given the
Normal increments $\Delta W^1(h)\coloneqq\Delta W^1(t_n,t_{n+1})$ and 
$\Delta W^2(h)\coloneqq\Delta W^2(t_n,t_{n+1})$ across $[t_n,t_{n+1}]$.

\subsection{\textbf{\textsf{Simulating Normal random variables}}}
We will start with the question: what is the most efficient method for generating
$\Delta W^1(h)$ and $\Delta W^2(h)$? The simple and direct answer is to
use the Matlab command 
\begin{quote}
\begin{center}
\texttt{sqrt(h)*randn}
\end{center}
\end{quote}
This command invokes an algorithm that has been scrupulously refined
and adapted over the years. One of the simplest efficient earlier incarnations 
of this algorithm is Marsaglia's polar method~\cite{Mar}.
The Box--M\"uller method is also very simple but not quite as efficient; 
see Kloeden and Platen~\cite{KP} for a discussion of these issues.
Also see Moro's inversion method~\cite{Moro}.
We outline Marsaglia's method here because of its simplicity and effectiveness. \medskip

\begin{algorithm}[\textbf{\textsf{Marsaglia's method}}] To produce two standard 
Normal samples:
\begin{enumerate}
\item Generate two independent uniform random samples
$U_1,U_2\in\text{\textsf{Unif}}([-1,1])$;
\item If $S\coloneqq U_1^2+U_2^2<1$ continue, otherwise repeat Step~1;
\item Compute $X_i=U_i/\sqrt{-2\,\ln(S)/S}$, for $i=1,2$; then $X_1$
and $X_2$ are independent standard Normal samples.
\end{enumerate}
\end{algorithm}

\subsection{\textbf{\textsf{Conditional distribution of L\'evy area}}}
The \emph{characteristic function} $\hat\phi$ of the probability density function for $\xi=A_{12}(h)$ 
given $\Delta W^1(h)$ and $\Delta W^2(h)$ is
\begin{equation*}
\hat\phi(\xi)=\frac{\tfrac12 h\xi}{\sinh(\tfrac12h\xi)}
\exp\Bigl(-\tfrac12a^2\bigl(\tfrac12h\xi\coth(\tfrac12h\xi)-1\bigr)\Bigr)
\end{equation*}
where $a^2=\bigl(\bigl(\Delta W^1(h)\bigr)^2+\bigl(\Delta W^1(h)\bigr)^2\bigr)/h$.
L\'evy derived this in a very succinct calculation in 1951; see L\'evy~\cite[pp.~171--3]{Levy}.
Since $\hat\phi$ is the characteristic function, i.e.\ the Fourier transform of the corresponding
probability density function, the actual probability density function $\phi$ is given
by the inverse Fourier transform (see for example Gaines and Lyons~\cite{GL:Mar}):
\begin{equation*}
\phi(x)=\tfrac{1}{\pi}\int_0^\infty \hat\phi(\xi)\,\cos(x\,\xi)\,\rd\xi.
\end{equation*}
The ungainly form of this probability density function means that generating
samples is not likely to be easy. For example, the 
simplest method for sampling from a continuous distribution $f$ is based
on the inversion of its (cumulative) distribution function 
$F(x)\coloneqq\int_{-\infty}^xf(\eta)\,\rd\eta$.
If we sample from the uniform distribution, say $U\sim\textsf{Unif}([0,1])$,
then $F^{-1}(U)$ is a sample from the target distribution. For this to be a practical
sampling method we must have an analytic form for $F$ or an
extremely efficient quadrature approximation for the integral in $F$ at our disposal. 
We don't have this for the probability density function of the L\'evy area $\phi$.

Several methods have been proposed for sampling from $\phi$. 
Gaines and Lyons~\cite{GL:Mar} proposed one of the most efficient,
based on Marsaglia's rectangle--wedge--tail method. However it can
be complicated to implement. Kloeden and Platen~\cite{KP} and 
Wiktorsson have proposed methods based on the Karhunen--Lo\`eve
expansion are much easier to code. Ryd\'en and Wiktorsson~\cite{RW}
proposed a method based on recognising the characteristic function 
$\hat\phi$ as a product of characteristic functions for a 
logistic random variable and an infinite sum of Poisson mixed Laplace
random variables. Gaines and Lyons~\cite{GL97} also
proposed a method based on the conditional expectation of the L\'evy area,
conditioned on intermediate Wiener increments.
We discuss these methods in the following four sections.
Stump and Hill~\cite{SH} have also proposed a very efficient
method, whose potential in a practical implementation is yet to be explored.

\begin{figure}
  \begin{center}
  \includegraphics[width=7.0cm,height=5.0cm]{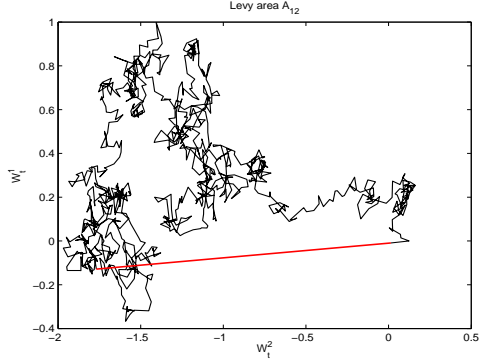}
  \end{center}
\caption{Sample two-dimensional Wiener path and enclosed chordal L\'evy area.}
\end{figure}

\subsection{\textbf{\textsf{Karhunen--Lo\`eve expansion method}}}
L\'evy~\cite{Levy} derived the form for the characteristic function $\phi$ 
for the L\'evy area, using the Karhunen--Lo\`eve expansion for a Brownian
bridge. This is an expansion in orthogonal basis functions. 
The details can be found in L\'evy~\cite{Levy} 
or Kloeden and Platen~\cite{KP}. If $U_{k},V_{k},X_{k},Y_{k}$ 
are independent $\textsf{N}(0,1)$ samples,
also independent of $\Delta W^1(h)$ and $\Delta W^2(h)$,
then the L\'evy area can be represented by
\begin{equation*}
A_{12}(h)=\frac{h}{2\pi}\sum_{k=1}^\infty 
\tfrac{1}{k}\Bigl(U_{k}\bigl(Y_{k}-\sqrt{\tfrac{2}{h}}\Delta W^2(h)\bigr)
-V_{k}\bigl(X_{k}-\sqrt{\tfrac{2}{h}}\Delta W^1(h)\bigr)\Bigr).
\end{equation*}
In practice, we truncate this expansion to only include $k\leqslant Q$ terms and
use the truncation, $\hat A_{12}(h)$, 
as an approximation for the L\'evy area. The important question
now, is as far as strong simulation is concerned, how many standard Normal
random variables do we need to simulate in order to have a sufficiently accurate
L\'evy area sample $\hat A_{12}(h)$, 
i.e.\/ how large must $Q$ be? Note that the coefficients in
the above expansion scale like $h/k$. 
The properties of the tail of the series, i.e.\/ all the terms for $k\geqslant Q+1$,
mean that it scales as $h/\sqrt{Q}$. For a Milstein numerical approximation we require
that the strong error is locally of order $h^{\frac{3}{2}}$; see Section~\ref{sec:strongerror}. 
Hence we must choose $Q\approx h^{-1}$ for a sufficiently accurate sample.

\subsection{\textbf{\textsf{Ryd\'en and Wiktorsson's method}}}
Ryd\'en and Wiktorsson~\cite{RW} proposed several methods, we detail here 
the most expedient.
The characteristic function $\hat\phi$ is the product of two characteristic functions
\begin{equation*}
\hat\phi_{X(h)}(\xi)=\frac{\tfrac12 h\xi}{\sinh(\tfrac12h\xi)}
\qquad\text{and}\qquad
\hat\phi_{Y(h)}(\xi)=\exp\Bigl(-\tfrac12a^2\bigl(\tfrac12h\xi\coth(\tfrac12h\xi)-1\bigr)\Bigr)
\end{equation*}
corresponding to the random variables $X(h)$ and $Y(h)$, respectively. We observe that
$\hat\phi_{X(h)}$ is the characteristic function of a logistic random variable which 
can be generated by the inverse method, i.e.\ pick $U\sim\textsf{Unif}([0,1])$ and 
let $X(h)=(h/2\pi)\,\log\bigl(U/(1-U)\bigr)$. Then using the identity 
\begin{equation*}
z\coth z-1=2\sum_{k=1}^\infty \frac{z^2}{\pi^2k^2+z^2},
\end{equation*}
we observe that 
\begin{equation*}
\hat\phi_{Y(h)}(\xi)=\exp\biggl(-a^2\sum_{k=1}^\infty\frac{\xi^2}{(2\pi k/h)^2+\xi^2}\biggr).
\end{equation*}
This can be viewed as a sum of compound Poisson random variables. Indeed if for each
$k\in\mathbb N$, we generate $N_k\sim\textsf{Poisson}(a^2)$ and then, 
for $j=1,\ldots,N_k$ generate independent Laplace random variables 
$Y_{jk}\sim\textsf{Laplace}(1/k)$, then
\begin{equation*}
Y(h)=\frac{h}{2\pi}\sum_{k=1}^\infty\sum_{j=1}^{N_k} Y_{jk},
\end{equation*}
has density $\phi_{Y(h)}$.
In a practical implementation we truncate this expansion to include $k\leqslant Q$
terms, and use the truncation as an approximation for the L\'evy area. 
Further the tail sum, by the central limit theorem, is asymptotically
Normally distributed and can be approximated by a Normal random variable.
This provides quite a dramatic improvement as 
it is possible to show that this method only requires the number
of standard Normal samples to be $Q\approx h^{-\frac{1}{2}}$.

\subsection{\textbf{\textsf{Wiktorsson's method}}}
Wiktorsson proposed a method that uses the Karhunen--Lo\`eve expansion method, 
but also simulates the tail sum as in the last method. 
Again, by the central limit theorem, the tail sum 
can be approximated by a Normal random variable, and the
corresponding improvement is that this method only requires the number
of standard Normal samples to be $Q\approx h^{-\frac{1}{2}}$.
Wiktorsson's method has been successfully implemented by 
Gilsing and Shardlow~\cite{GilSh} in their SDELab, to where
the interested reader is referred.

\subsection{\textbf{\textsf{Conditional expectation}}}
One more approach to simulating the L\'evy area,
or equivalently $J_{12}(t_n,t_{n+1})$, is
based on replacing $J_{12}(t_n,t_{n+1})$, by its conditional expectation $\hat J_{12}(t_n,t_{n+1})$,
as follows. Suppose we are about to perform the numerical update
for the solution across the interval $[t_n,t_{n+1}]$.
We generate $Q$ pairs of independent standard Normal random variables 
$X_q,Y_q\sim \textsf{N}(0,1)$ for $q=1,\ldots,Q$. 
Set $\tau_q=t_n+q\Delta t$, for $q=0,\ldots,Q-1$, and 
$\Delta W^1(\tau_q)=\sqrt{\Delta t}\,X_q$ and $\Delta W^2(\tau_q)=\sqrt{\Delta t}\,Y_q$,
where $\Delta t$ is defined by $Q\Delta t=h$. We thus generate a two-dimensional
Brownian sample path on $[t_n,t_{n+1}]$. We can take 
$\Delta W^1(h)$ and $\Delta W^2(h)$ to be
the increments across the interval $[t_n,t_{n+1}]$. More importantly,
we can use the intervening path information we have generated 
on the scale $\Delta t$ to approximate $J_{12}(t_n,t_{n+1})$. 
Indeed, $J_{12}(t_n,t_{n+1})$ can be expressed as 
\begin{align*}
J_{12}(t_n,t_{n+1})=&\; \int_{t_n}^{t_{n+1}}\int_{t_n}^\tau \,\mathrm{d}W_{\tau_1}^1
\,\mathrm{d}W_{\tau}^2\\
=&\; \sum_{q=0}^{Q-1}\int_{\tau_q}^{\tau_{q+1}} (W_{\tau}^1-W_{\tau_q}^1)
+(W_{\tau_q}^1-W_{t_n}^1)
\,\mathrm{d}W_{\tau}^2\\
=&\;\sum_{q=0}^{Q-1}J_{12}(\tau_q,\tau_{q+1})
+\sum_{q=0}^{Q-1}\bigl(W_{\tau_q}^1-W_{t_n}^1\bigr)
\,\Delta W^2(\tau_q).
\end{align*}
The quantity 
\begin{equation*}
\hat J_{12}(t_n,t_{n+1})\coloneqq\sum_{q=0}^{Q-1}\bigl(W_{\tau_q}^1-W_{t_n}^1\bigr)
\,\Delta W^2(\tau_q)
\end{equation*}
represents the expectation of $J_{12}(t_n,t_{n+1})$ conditioned on the 
increments $\Delta W^1(\tau_q)$ and $\Delta W^2(\tau_q)$.
From an algebraic and geometric perspective, $\hat J_{12}(t_n,t_{n+1})$ 
represents a suitable approximation to $J_{12}(t_n,t_{n+1})$.
Computing its mean-square strong error we see that
\begin{equation*}
\left\|J_{12}(t_n,t_{n+1})-\hat J_{12}(t_n,t_{n+1})\right\|_{L_2}^2
=\sum_{q=0}^{Q-1}\bigl\|J_{12}(\tau_q,\tau_{q+1})\bigr\|_{L^2}^2
=Q(\Delta t)^2=h^2/Q.
\end{equation*}
Hence its root-mean-square strong error is $h/\sqrt{Q}$.
Thus, as for the Karhunen--Lo\`eve expansion approach, 
to achieve a suitable approximate sample for the stochastic area integral, 
this method requires $Q\approx h^{-1}$.
One advantage of this method is that it is very convenient for generating
log-log error plots. 
\medskip


\section{\textsf{Strong error}}\label{sec:strongerror}
We will focus here on the global, strong $L^2$ error. In practical terms,
the global error is generated by the accumulation of contributions from
the local error. The local error is itself the leading order terms 
in the remainder $R_t$, say, of our truncated stochastic Taylor series.
Note that there is also a contribution to the global error from the
approximate simulation of the L\'evy area, however we will assume here,
that the L\'evy area has been sufficiently accurately simulated, as discussed
in detail in the last section, so that its contribution is \emph{small},
in comparison to the truncation error.

Suppose we base a strong numerical approximation on truncating the stochastic
Taylor expansion (in Stratonovich form). Let $\hat y$ denoted the
truncated expansion and $R$ the corresponding remainder; hence the 
exact solution is $y=\hat y+R$. To guarantee our numerical scheme
based on such a truncation is globally of order $h^{m}$, where $m\in\mathbb Z/2$,
which terms must we keep in $\hat y$? We give the following rule.\medskip

\noindent\textsf{Rule of Thumb:} Terms in the remainder $R$ of $L^2$ measure
$h^m$ with:
\begin{itemize}
\item zero expectation, accumulate so as to contribute to the global
error as $h^{m-\frac12}$ order terms;
\item non-zero expectation, accumulate so as to contribute to the global
error as $h^{m-1}$ order terms.
\end{itemize}

Hence to achieve an integrator with global error of order $h^{m}$,
we must retain in $\hat y$:
\begin{itemize}
\item all terms with $L^2$ measure of order $h^{m'}$ for all $m'\leqslant m$;
\item the expectation of all terms of order $h^{m+\frac12}$ which have non-zero
expectation (the corresponding terms left in remainder will then have
zero expectation).
\end{itemize}\medskip

\noindent\textsf{Example (Euler--Maruyama).} Recall that we based
the Euler--Maruyama approximation on the truncated It\^o Taylor
series. If we had truncated the stochastic Taylor series in
Stratonovich form, then according to the rules above, to achieve global order $h^{\frac12}$
we should retain in our integrator the expectation of the terms 
\begin{equation*}
\sum_{i=1}^d(V_i\cdot\partial_yV_i)\,J_{ii}.
\end{equation*}
Since $\Eb\, J_{ii}=\tfrac12 h$, we thus recover 
the corresponding truncated It\^o Taylor series.\smallskip

\section{\textsf{Further issues}}
There are many important simulation issues we have not had space
to discuss. Chief among these is the numerical stability of the 
strong methods we have explicitly outlined. This issue is discussed
in Higham~\cite{Higham:intro} and more can be found for example in 
Buckwar, Horv\'ath--Bokor and Winkler~\cite{BHBW}.


\appendix

\section{\textsf{Stratonovich to It\^o relations}}
We give here some Stratonovich to It\^o relations for 
convenience for the reader---more 
details can be found in Kloeden and Platen~\cite{KP}. For the 
words $w$ shown, the Stratonovich integrals $J_w$ can be 
expressed in terms of It\^o integrals $I_w$ as follows:
\begin{align*}
w=a_1a_2\colon&\; J_w=I_w+\tfrac12I_0\,\delta_{a_1=a_2\neq0};\\
w=a_1a_2a_3\colon&\; J_w=I_w+\tfrac12(I_{0a_3}\,\delta_{a_1=a_2\neq0}
+I_{a_10}\,\delta_{a_2=a_3\neq0});\\
w=a_1a_2a_3a_4\colon&\; J_w=I_w
+\tfrac14 I_{00}\,\delta_{a_1=a_2\neq0}\delta_{a_3=a_4\neq0}\\
&\;\qquad\qquad+\tfrac12(I_{0a_3a_4}\,\delta_{a_1=a_2\neq0}+I_{a_10a_4}\,\delta_{a_2=a_3\neq0}
+I_{a_1a_20}\,\delta_{a_3=a_4\neq0}).
\end{align*}
Note that the expectation of any It\^o integral $I_w$ is zero, i.e.\/: 
$\Eb\, I_w=0$ for any word $w\in\Ab^\ast$ which has at least one non-zero 
letter.

\section{\textsf{Sample program for weak and strong Euler--Maruyama}}

We provide the listing for the weak vs strong Euler--Maruyama simulation shown
in Figure~\ref{fig:weakstrongsim}.

\begin{lstlisting}[frame=topline,caption={Weak vs strong simulation},label=WVSS]

%%    Weak and strong simulation example plot
%%    Euler-Maruyama approximations
%%    Example scalar linear SDE

%%    Parameter values

a=3.0;                              % drift coefficient
b=1.4;                              % diffusion coefficient is b
y0=1.0;                             % initial data

P=10;                               % total # of sample paths
h=0.05;                             % stepsize
T=1.0;                              % global time interval
N=T/h;                              % number of subintervals

%%    Binomial branching process increments

dw=zeros(N,P);
w=zeros(N+1,P);
binom=binornd(1,1/2,[N,P]);         % Gives 0 or 1 with prob 1/2
dw=sqrt(h)*(1-2*binom);             % Binomial increments dw
w(2:N+1,:)=cumsum(dw,1);            % Binomial paths themselves

%%    Weak solution by Euler-Maruyama approximation

y=zeros(N+1,P);                     % y is weak solution

for p=1:P                           % set loop for each path
y(1,p)=y0;                          % initial data

for n=1:N
    y(n+1,p)=y(n,p)+a*y(n,p)*h+b*y(n,p)*dw(n,p);    
end

end

%%    Approximate Wiener path increments

dW=zeros(N,P);
W=zeros(N+1,P);
dW=sqrt(h)*randn(N,P);              % Brownian increments dW
W(2:N+1,:)=cumsum(dW,1);            % Brownian paths themselves

%%     Strong solution by Euler-Maruyama approximation

Y=zeros(N+1,P);                     % Y is strong solution

for p=1:P                           % set loop for each path
Y(1,p)=y0;                          % initial data

for n=1:N
    Y(n+1,p)=Y(n,p)+a*Y(n,p)*h+b*Y(n,p)*dW(n,p);    
end
end

%%     Compute the expectations of y and Y at each timestep

expect_y=zeros(N+1,1);
expect_Y=zeros(N+1,1);

for n=1:N+1
    expect_y(n)=mean(y(n,:));
    expect_Y(n)=mean(Y(n,:));
end
\end{lstlisting}

\section{\textsf{Example strong simulation program}}

\subsection{\textbf{\textsf{Heston model strong simulation}}}
\begin{figure}
  \begin{center}
  \includegraphics[width=12.0cm,height=5.0cm]{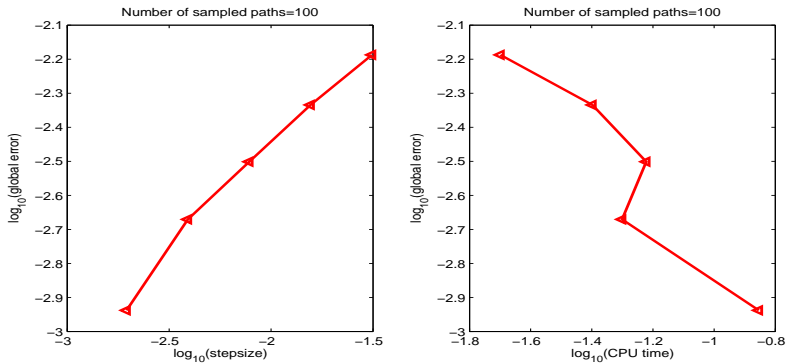}
  \end{center}
\caption{Error vs stepsize and error vs CPU time for the Heston model.
The parameter values can be seen in the program listing.}
\label{fig:simsde0p5}
\end{figure}

We provide here a sample program that shows how to perform log-log
error plots for a strong simulation. We used a real example,
the Heston model, and applied the full truncation Euler--Maruyama
type numerical scheme devised by Lord, Koekkoek and Van Dijk~\cite{LKVD}.
The log-log error vs stepsize, and error vs CPU time, are shown
in Fig.~\ref{fig:simsde0p5}. Note that to estimate the strong 
global error, we must compare solutions for different stepsizes
along the \emph{same} path, before taking the expectation.

\begin{lstlisting}[frame=topline,caption={Heston model strong simulation},label=HMSS]

%%       Heston model integration       

% Global data

T0=0;                            % integration interval...
T=1;                             % ...is [T0,T]
M=10;
hmin=(T-T0)/2^M;                 % smallest timestep
Mstart=4; 
hmax=(T-T0)/2^Mstart;            % largest timestep 
Q=(1/(hmin))*(hmax/hmin);        % # of quad pts for hmin
dt=hmax/Q;                       % quad scale for hmin
R=M-Mstart+1;                    % # of solution approximations...
                                 % ...at different stepsizes
P=100;                           % total # of Brownian paths
alpha=2.0;                       % parameters
theta=0.09;
beta=0.1;
rho=0.5;
mu=0.05;
ic=[1.0; 0.09];                  % initial data

YFT=zeros(P,R,2);                % approximate solution
clockYFT=zeros(1,R);             % CPU timings

for p=1:P
    for r=1:R
        YFT(p,r,:)=ic;           % make sure start with IC
    end
end

for jj=1:2^Mstart
YFTold=YFT;

for p=1:P                        % loop for each path

%%       Start loop computing over intervals of length hmax

siv=hmax/hmin;

%%       Generate dW1 and dW2 on smallest scale

    dW1=sqrt(dt)*randn(1,Q);     % Brownian increments dW1
    dW2=sqrt(dt)*randn(1,Q);     % Brownian increments dW2
    dW0=(zeros(1,Q)+1)*dt;

%%       Start loop computing Js with different stepsizes     

    for r=1:R
    SF=2^(r-1);                  % scale factor for stepsizes...
    h=SF*hmin;                   % stepsize
    L=hmax/h;                    % # of timesteps needed
    QR=Q/(SF^2);                 % # of quadrature steps (total)...

% Compute dw0, dw1, dw2, w0, w1, w2   

     dw0=zeros(1,QR);
     dw1=zeros(1,QR);
     dw2=zeros(1,QR);
     w0=zeros(1,QR);
     w1=zeros(1,QR);
     w2=zeros(1,QR);

     for j=1:QR
         dw0(j)=sum(dW0((j-1)*(SF^2)+1:j*(SF^2)));
         dw1(j)=sum(dW1((j-1)*(SF^2)+1:j*(SF^2)));
         dw2(j)=sum(dW2((j-1)*(SF^2)+1:j*(SF^2)));
     end

     QF=QR/L;                    % # of quadrature steps in h...

     for n=1:L
         w0((n-1)*QF+1:n*QF)=cumsum(dw0((n-1)*QF+1:n*QF));
         w1((n-1)*QF+1:n*QF)=cumsum(dw1((n-1)*QF+1:n*QF));
         w2((n-1)*QF+1:n*QF)=cumsum(dw2((n-1)*QF+1:n*QF));
     end
   
    oldclockYFT=clockYFT(r);
    ts=cputime;
    YFT(p,r,:)=YFTapprox(p,h,L,T0,QF,dw0,dw1,dw2,w0,w1,w2, ...
                                 alpha,theta,beta,rho,mu,...
                                 YFTold(p,r,:));
    clockYFT(r)=cputime-ts+oldclockYFT;

    end
end
end

stepsizes=log10((2.^([1:R-1]))*hmin);

save('stepsizes','stepsizes')
save('P','P')
save('YFT','YFT')
save('clockYFT','clockYFT')
\end{lstlisting}

\subsection{\textbf{\textsf{Program listing: integrator}}}
We give the program listing for the Heston model full truncation
Euler--Maruyama integrator.

\begin{lstlisting}[frame=topline,caption={Full truncation
Euler--Maruyama integrator},label=FTEM]

%%    Full truncation for volatility and exponential for index: 

function trunc=YFTapprox(p,h,L,T0,QR,dW0,dW1,dW2,...
                         W0,W1,W2,alpha,theta,beta,rho,mu,ic)

trunc=zeros(2,1);

% Compute dJ1, dJ2 on scale h=SF*hmin

J1=zeros(1,L);
J2=zeros(1,L);
J1(1)=W1(QR);
J2(1)=W2(QR);
pts=2*QR:QR:L*QR;
J1(2:L)=W1(pts);
J2(2:L)=W2(pts);

%%        initial data

S=ic(1);                                  % asset price
v=ic(2);                                  % volatility

for n=1:L
    Sold=S;
    vold=v;
    S=exp((mu-max(0,vold)/2)*h+sqrt(max(0,vold))*J1(n))*Sold;
    v=vold+alpha*(theta-max(0,vold))*h+beta*(rho*J1(n)...
          +sqrt(1-rho^2)*J2(n))*sqrt(max(0,vold));
end

trunc=[S; v];

end
\end{lstlisting}

\subsection{\textbf{\textsf{Program listing: log-log strong error plots}}}
The following program performs the log-log plots for the strong $L^2$ error measure.

\begin{lstlisting}[frame=topline,caption={Log-log strong error plots},label=LLEP]

%%     Log-log strong error plots

load stepsizes
load YFT
load clockYFT
load P

R=length(stepsizes);
errorYFT=zeros(1,R);
diffYFT=zeros(P,R);

% the "norm" below is the Euclidean norm

for r=1:R
    for p=1:P
    diffYFT(p,r)=norm(YFT(p,r+1)-YFT(p,1));
    end
end

% now take the L^2 norm measure 

errorYFT=sqrt(mean(diffYFT.^2,1));

figure
subplot(1,2,1)
plot(stepsizes(1:end-1),log10(errorYFT(1:end-1)),...
                '-ks','LineWidth',2)
xlabel('log_{10}(stepsize)')
ylabel('log_{10}(global error)')
title(['Number of sampled paths=',int2str(P)])
subplot(1,2,2)
plot(log10(clockYFT(2:end-1)),log10(errorYFT(1:end-1)),...
                 '-ks','LineWidth',2)
xlabel('log_{10}(CPU time)')
ylabel('log_{10}(global error)')
title(['Number of sampled paths=',int2str(P)])
\end{lstlisting}

\end{document}